\documentclass[12pt,a4paper]{article}
\usepackage{epsfig,latexsym,amsfonts,amssymb,amsmath,amscd,
graphics,theorem} %\usepackage[notref,notcite]{showkeys}
\theoremstyle{plain}
\newtheorem{lemma}{Lemma}%[section]

\newtheorem{theorem}[lemma]{Theorem}\newtheorem{conjecture}[lemma]{Conjecture}

{\theorembodyfont{\rmfamily} 
\font\ncsc=cmcsc10  \font\ntt=cmtt12
\begin{document}
\baselineskip=15pt
\newcommand{\pperp}{\hbox{$\perp\hskip-6pt\perp$}}
\newcommand{\ssim}{\hbox{$\hskip-2pt\sim$}}
\newcommand{\N}{{\mathbb N}}\newcommand{\Delp}{{\Pi}}
\newcommand{\A}{{\mathbb A}}
\newcommand{\Z}{{\mathbb Z}}
\newcommand{\R}{{\mathbb R}}
\newcommand{\C}{{\mathbb C}}
\newcommand{\Q}{{\mathbb Q}}
\newcommand{\PP}{{\mathbb P}}
\newcommand{\mnote}{\marginpar}
\newcommand{\Id}{{\operatorname{Id}}}
\newcommand{\oeps}{{\overline\eps}}
\newcommand{\oDel}{{\widetilde\Del}}
\newcommand{\real}{{\operatorname{Re}}}
\newcommand{\conv}{{\operatorname{conv}}}
\newcommand{\Span}{{\operatorname{Span}}}
\newcommand{\Ker}{{\operatorname{Ker}}}
\newcommand{\Hyp}{{\operatorname{Hyp}}}
\newcommand{\Fix}{{\operatorname{Fix}}}
\newcommand{\sign}{{\operatorname{sign}}}
\newcommand{\Tors}{{\operatorname{Tors}}}
\newcommand{\oi}{{\overline i}}
\newcommand{\oj}{{\overline j}}
\newcommand{\ob}{{\overline b}}
\newcommand{\os}{{\overline s}}
\newcommand{\oa}{{\overline a}}
\newcommand{\oy}{{\overline y}}
\newcommand{\ow}{{\overline w}}
\newcommand{\ou}{{\overline u}}
\newcommand{\ot}{{\overline t}}
\newcommand{\oz}{{\overline z}}
\newcommand{\bw}{{\boldsymbol w}}
\newcommand{\bx}{{\boldsymbol x}}
\newcommand{\bu}{{\boldsymbol u}}
\newcommand{\bz}{{\boldsymbol z}}
\newcommand{\eps}{{\varepsilon}}
\newcommand{\proofend}{\hfill$\Box$\bigskip}
\newcommand{\Int}{{\operatorname{Int}}}
\newcommand{\pr}{{\operatorname{Pr}}}
\newcommand{\grad}{{\operatorname{grad}}}
\newcommand{\rk}{{\operatorname{rk}}}
\newcommand{\im}{{\operatorname{Im}}}
\newcommand{\sk}{{\operatorname{sk}}}
\newcommand{\const}{{\operatorname{const}}}
\newcommand{\Sing}{{\operatorname{Sing}}}
\newcommand{\conj}{{\operatorname{Conj}}}
\newcommand{\Pic}{{\operatorname{Pic}}}
\newcommand{\Crit}{{\operatorname{Crit}}}
\newcommand{\Ch}{{\operatorname{Ch}}}
\newcommand{\discr}{{\operatorname{discr}}}
\newcommand{\Tor}{{\operatorname{Tor}}}
\newcommand{\Conj}{{\operatorname{Conj}}}
\newcommand{\val}{{\operatorname{val}}}
\newcommand{\Val}{{\operatorname{Val}}}
\newcommand{\defect}{{\operatorname{def}}}
\newcommand{\tmu}{{\C\mu}}
\newcommand{\ov}{{\overline v}}
\newcommand{\ox}{{\overline{x}}}
\newcommand{\tet}{{\theta}}
\newcommand{\Del}{{\Delta}}
\newcommand{\bet}{{\beta}}
\newcommand{\kap}{{\kappa}}
\newcommand{\del}{{\delta}}
\newcommand{\sig}{{\sigma}}
\newcommand{\alp}{{\alpha}}
\newcommand{\Sig}{{\Sigma}}
\newcommand{\Gam}{{\Gamma}}
\newcommand{\gam}{{\gamma}}
\newcommand{\Lam}{{\Lambda}}
\newcommand{\lam}{{\lambda}}
\newcommand{\SC}{{SC}}
\newcommand{\MC}{{MC}}
\newcommand{\nek}{{,...,}}
\newcommand{\cim}{{c_{\mbox{\rm im}}}}
\newcommand{\mathto}{\mathop{\to}}
\newcommand{\op}{{\overline p}}

\newcommand{\w}{{\omega}}

\title{Logarithmic equivalence
of Welschinger and Gromov-Witten invariants}
\author{Ilia Itenberg \and Viatcheslav Kharlamov
\and Eugenii Shustin}
\date{}
\maketitle

\centerline{\it\small To the memory of Andrey Bolibruch, a lively
man of creative mind and open soul}

\vskip1cm

\begin{abstract}
The
Welschinger numbers, a kind of a real analog of the
Gromov-Witten numbers which count the complex rational curves
through a given generic collection of points,
bound from below the number of real rational curves for any real
generic collection of points.
By the
logarithmic equivalence of sequences we mean the asymptotic equivalence
of their logarithms.
We prove such an equivalence for
the Welschinger and Gromov-Witten numbers of
any toric Del Pezzo surface with
its tautological real structure, in particular,
of the projective
plane, under the hypothesis that all, or almost all,
chosen points are real.
We also study the positivity
of Welschinger numbers
and their monotonicity
with respect to the number of
imaginary points.
\end{abstract}

\section{Introduction}\label{intro}

The present note is devoted to
an asymptotic enumeration of real rational curves interpolating
fixed real collections of points
in a real surface $\Sigma$,
more precisely, to the following
question: {\it given a real divisor $D$
and a generic collection $\bw$ of $c_1(\Sigma)\cdot D-1$ real points in $\Sigma$,
how many
of the complex
rational curves in the linear system $|D|$
passing through $\bw$
are real \/}? By rational curves we mean irreducible genus zero
curves and
their degenerations,
so that they form in $|D|$ a projective subvariety $S(\Sigma, D)$;
this subvariety is called
the
{\it Severi variety}.
A curve on a real surface
$\Sig$
is called real, if it is invariant under the involution
$c:\Sigma\to\Sigma$ defining
the real structure of $\Sigma$.

While, under mild conditions on $\Sigma$ and $D$, the number of
complex curves in question is the same for all generic $\bw$ (it
equals to the degree of $S(\Sigma, D)$), it is no more the case
for real curves (except few very particular situations). For
example, if $\Sigma$ is the projective plane, such a
non-invariance manifests starting from the degree $d=3$: if $d=3$
and all the $3d-1=8$ points are real, the number of interpolating
real rational curves takes the values $8, 10$, and $12$ (twelve is
the number of complex interpolating curves).

Till recently, in this interpolation problem, even the existence
of at least one real solution for arbitrary choice of a generic
collection of $3d-1$ points in $\R P^2$ was known only for $d\le
3$. The situation has radically changed after the discovery by
J.-Y.~Welschinger~\cite{W, W1} of a way to attribute weights of
$\pm 1$ to real solutions so that the number of solutions counted
with weights is independent of the configuration of points; more
precisely, a real collection $\bw$ may even contain complex
conjugated pairs of points, and then the result depends only on
the number $m$ of pairs of imaginary points in $\bw$, and $d$. As
an immediate consequence, the absolute value of Welschinger's
invariant $W_{d,m}$ provides a lower bound on the number
$R_{d,m}(\bw)$ of real solutions: $R_{d,m}(\bw) \ge |W_{d,m}|$.

In \cite{IKS} we proved
an inequality $W_{d, 0}\ge \frac12 d!$
which, due to Welschinger's
lower bound,
implies that
{\it for any
integer $d\ge 1$, through any $3d-1$ generic points in $\R P^2$
there can be traced at least $\frac 12 d!$ real rational curves of degree
$d$.}

Comparing this lower bound with the number $N_d$ of complex plane
rational curves of degree~$d$ passing through $3d - 1$ generic
points one can observe that the logarithm of the bound is
asymptotically equal to $\frac{1}{3}\log N_d$. In fact, $\log
d!\sim d\log d $, and $\log N_d\sim 3d\log d$, as follows from the
inequalities $(3d-4)!\cdot 54^{-d}\le N_d \le (3d-5)!$ which in
their turn follow from Kontsevich's recurrent formula
{\rm\cite{KM}} (a more precise asymptotics is found in \cite{Itz},
cf. Remark (3) in \ref{statements} below). The next natural
question arises: {\it how far are $W_d = W_{d, 0}$ from $N_d$ in
the logarithmic scale\rm }?

In
the present
note
we give a complete answer to the above question.

\begin{theorem}\label{t3}
$\;$ The sequences $\log W_d$ and $\log N_d$ are asymptotically
equivalent.
\end{theorem}

To establish this asymptotic equivalence
we improve the lower bound
$W_{d, 0}\ge \frac12 d!$ we
obtained in
\cite{IKS}.
As in \cite{IKS},
the method we use is based on Mikhalkin's
approach~\cite{Mi}
to counting nodal curves passing through specific
configurations of points,
an approach which
deals with a corresponding count of
tropical curves.
In fact,
we
slightly modify the construction of~\cite{IKS} in order to
multiply the
logarithm of the
bound by $3$.

We
obtain similar
results
for
curves on
other toric Del Pezzo surfaces
equipped with
their
tautological
real structure.
We also study the behavior of the double sequence $W_{d,m}$ with
respect to $m$. To treat $W_{d,m}$ with $m>0$ we use Shustin's
counting scheme, see \cite{Sh1}, which extends Mikhalkin's scheme
from pure real data ($m=0$) to arbitrary complex conjugation
invariant ones. (It
may be worth noticing, that these counting schemes can be
considered as an advanced, almost explicit, version of the Viro
patchworking, cf. \cite{Mi, Sh0, Sh1, Vi}.)

We also note that our asymptotic results give some information on
the convergency domain of the Gromov-Witten potential, see Remarks
in \ref{statements}.

The paper is organized as follows. A separate section,
Section~\ref{sec2}, is devoted to the case where all the points
are real. The main results are summarized in Theorem \ref{t5}
(which contains the above Theorem \ref{t3} as a particular case).
The proof is divided in two parts: an upper bound for
Gromov-Witten invariants (Lemma \ref{upperbound})
and a lower
bound for Welschinger invariants (Theorem \ref{kb}). In
Section~\ref{further} we analyze the case where some of the points
are imaginary. We start from explicit calculations of Welschinger
invariants in few particular cases. Then, after resuming Shustin's
general counting scheme, we derive from it few results on
positivity, monotonicity, and asymptotics of Welschinger
invariants, Theorems \ref{tn4} and \ref{tn5}.

As is already mentioned,
the counting scheme used in Section 2 is the
same as in \cite{IKS}. It is taken from \cite{Mi0,Sh0} and its
summary
is found in \cite{IKS}.
At the same time,
it is included as a
special case into Shustin's extended scheme described in
Section 3, so that a reader can
reconstruct the first
scheme from the extended one
specializing
the number of imaginary points
to zero.

\bigskip

{\bf Acknowledgements}. The first and the second authors are
members of Research Training Network RAAG CT-2001-00271; the
second author is as well a member of Research Training Network
EDGE CT-2000-00101.
The main results of this note were obtained
during the stay of the authors at the Mathematical Sciences
Research Institute in Berkeley, and we thank this institution for
the hospitality and excellent work conditions.

\section{Asymptotics of the Welschinger invariants
for purely real data}\label{sec2}

\hfill
\begin{picture}(50,25)
\put(0,0){\line(1,0){50}}
\put(50,0){\line(0,1){25}}
\put(0,25){\line(1,0){50}}
\put(0,0){\line(0,1){25}}
\end{picture}

{\bf\tiny
\hfill He had bought a large map representing the sea,
\vskip-.2cm
\hfill Without the least vestige of land:
\vskip-.2cm
\hfill And the crew were much pleased when they found it to be
\vskip-.2cm
\hfill A map they could all understand.
\vskip-.2cm}

{\sc\tiny \hfill The Hunting of the Snark, Lewis Carroll}

\subsection{Notations}\label{notations}
There are five \emph{unnodal}
({\it i.e.},
not
containing
any $(-2)$-curve)
toric Del Pezzo surfaces: the projective plane
$\PP^2$, the product
of projective lines
$Q=\PP^1 \times \PP^1$, and $\PP^2$ with
$k$ blown up generic points, $k=1,2$ or 3; the latter three surfaces
are denoted by $P_k$.
Let $E_1, \ldots ,E_k$ be the
exceptional divisors of $P_k\to \PP^2$ and $L\subset P_k$
the pull back of a generic straight line.

We equip $\PP^2$ and $Q$ with their tautological real and toric
structures. For $P_k$, we choose the blown up points in $\PP^2$ to
be among the three $0$-dimensional orbits, so that $P_k$ inherits
toric and real structures from $\PP^2$.

Let $D$ be an ample divisor on $\Sigma$. The linear system $|D|$
is generated, with respect to suitable real toric coordinates, by
monomials $x^iy^j$, where $(i,j)$ ranges over all the integer
points ({\it i.e.}, points having integer coordinates) of a convex
polygon~$\Delp=\Delp_D$ of the following form. If $\Sig = \PP^2$
and $D = d[\PP^1]$, then $\Delp$ is the triangle with vertices
$(0,0)$, $(d,0)$, and $(0,d)$. If $\Sig=Q$ and $D$ is of bi-degree
$(d_1, d_2)$, then $\Delp$ is the parallelogram with vertices
$(d_2,0)$, $(d_1+d_2,0)$, $(d_1,d_2)$, and $(0,d_2)$. If
$\Sig=P_k$, $k=1,2,3$, and $D = dL-\sum_{i=1}^k d_iE_i$, then
$\Delp$ is respectively the trapeze with vertices $(0,0)$,
$(d,0)$, $(d_1, d-d_1)$, $(0,d-d_1)$, or the pentagon with
vertices $(0,0)$, $(d-d_1,0)$, $(d-d_1,d_1)$, $(d_2,d-d_2)$,
$(0,d-d_2)$, or the hexagon with vertices $(d_3,0)$, $(d-d_1,0)$,
$(d-d_1,d_1)$, $(d_2,d-d_2)$, $(0,d-d_2)$, $(0,d_3)$. (A choice of
a parallelogram instead of the rectangle in the case $\Sig=Q$
simplifies the analysis of the irreducibility of the curves
appearing in the proof of Theorem \ref{kb}.)

Let $r = r(\Delp)$ be the number of integer points on the boundary
$\partial \Delp$ of $\Delp$ diminished by $1$, and $\delta(\Delp)$
be the number of interior integer points of~$\Delp$. Note that
$r(\Delp) = c_1(\Sig)\cdot D-1$ and $\delta(\Delp)$ is the genus
of nonsingular representatives of $|D|$. As is well known, the
number of curves of genus $0\le g\le\delta(\Delp)$ in $|D|$
passing through $c_1(\Sig)\cdot D-1+g$ generic points is finite.
Denote by $N_{\Sig, D}$, or shortly $N_D$, the number of complex
rational curves in $|D|$ passing through $r$ given generic points
of $\Sigma$ (note that rational curves in $|D|$ which pass through
$r$ generic points are irreducible and nodal). Due to the
genericity of the complex structure of unnodal Del Pezzo surfaces,
these enumerative numbers coincide with the Gromov-Witten genus
zero invariants. An inductive procedure reconstructing their
values was given for $\Sig=\PP^2$ by M.~Kontsevich and for other
unnodal Del Pezzo surfaces by M.~Kontsevich and Yu.~Manin, see
\cite{KM}.

\subsection{Welschinger numbers}\label{Welsch-inv}

Let us fix an integer $m$ such that $0\le 2m\le r$ and introduce a
real structure $c_{r,m}$ on $\Sig^r$ which maps $(z_1, \ldots ,
z_r)\in\Sig^r$ to $(z'_1, \ldots , z'_r)\in\Sig^r$ with $z'_i =
c(z_i)$ if $i>2m$, and $(z'_{2j-1}, z'_{2j}) = (c(z_{2j}),
c(z_{2j-1}))$ if $j\le m$. With respect to this real structure, a
point $\bw=(z_1,...,z_r)$ is real, {\it i.e.},
$c_{r,m}$-invariant, if and only if $z_i$ belongs to the real part
$\R\Sig$ of $\Sig$ for $i>2m$ and $z_{2j-1},z_{2j}$ are conjugate
for $j \le m$. In what follows we work with
the
open dense subset
$\Omega_{r,m}(\Sig)$ of $\R\Sig^r=\Fix\, c_{r,m}$ constituted of
$c_{r,m}$-invariant $r$-tuples $\bw=(z_1,...,z_r)$ with pairwise
distinct $z_i\in\Sig$.

By abuse of language, we say that a curve $C$ in $\Sig$ passes
through $\bw\in\Sig^r$ if $C$ contains all the components
$z_i\in\Sig$ of $\bw$.

In the spirit of \cite{W,W1}, given a {\bf generic}
$\bw\in\Omega_{r,m}(\Sig)$, we introduce the number
$W_{D,m}^{even}(\bw)$ (resp., $W_{D,m}^{odd}(\bw)$) of irreducible
real rational curves in $|D|$ passing through $\bw$ and having
even (resp., odd) number of solitary nodes ({\it i.e.}, real
double points, where a local equation of the curve can be written
over $\R$ in the form $x^2+y^2=0$).  The {\it Welschinger number}
$W_{D,m}(\bw)$ is defined by $W_{D,m}(\bw)= W_{D,m}^{even}(\bw)-
W_{D,m}^{odd}(\bw)$.\footnote{In \cite{IKS} we considered only the
case $m=0$ and used a slightly different notation.}

\begin{theorem}\label{p1}
{\rm (J.-Y.~Welschinger, see~\cite{W, W1})}. The value
$W_{D,m}(\bw)$ does not depend on the choice of a {\rm
}generic{\rm } element $\bw$ in $\Omega_{r,m}(\Sig)$.
\end{theorem}

In the case $m = 0$, the number $W_{D,m} = W_{D,m}(\bw)$ is
denoted by $W_D$. We also use more detailed notations
$W_{\Sigma,D}$ and $W_{\Sigma,D,m}$ when we need to work with
several surfaces $\Sigma$ simultaneously.

\subsection{Key bound}\label{key}
The following bound plays a crucial
role in our treatment of logarithmic equivalence of Welschinger
and Gromov-Witten invariants.

\begin{theorem}\label{kb}
Let  $\Sig$ be
$\PP^2$,
$Q$,
$P_1$, $P_2$, or $P_3$
equipped with its tautological real structure,
and
$D$
an ample divisor on $\Sigma$. Then
\begin{equation}\label{kbound}
W_{nD}\ge \exp (a\, n\log{n}+O(n)),\quad n\in\N, \quad a=c_1(\Sig)\cdot D\ .
\end{equation}
\end{theorem}

\subsection{Main corollary}
\label{statements}

Let as above $\Sig$ be $\PP^2$, $Q$, $P_1$, $P_2$, or $P_3$
equipped with its tautological real structure, and $D$ an ample
divisor on $\Sigma$. We prove the following theorem.

\begin{theorem}\label{t5}
The sequences $\log W_{nD}$ and $\log N_{nD}$, $n\in\N$, are
asymptotically equivalent. More precisely,
\begin{equation}\label{main}
\log W_{nD} = \log N_{nD} + O(n)\quad\text{and}\quad \log N_{nD} =
(c_1(\Sigma) \cdot D)\cdot n\log n + O(n).
\end{equation}
\end{theorem}

{\bf Proof}. Due to
$N_{nD}\ge W_{nD}$, all the statements
follow from Theorem \ref{kb} and Lemma \ref{upperbound}.
\proofend

\begin{lemma}\label{upperbound} The following inequality holds:
\begin{equation}\label{upperb}
\log N_{nD}\le (c_1(\Sigma) \cdot D)\cdot n\log n + O(n).
\end{equation}
\end{lemma}

{\bf Proof}.
The case of $\Sig=Q$ with a
divisor $D$ of bi-degree $(d_1,d_2)$ is equivalent to the
case of $\Sig=P_2$  with $D= dL-d_1E_1-d_2E_2$,
$d=d_1+d_2$.
The case of $\Sig=\PP^2$ is covered by \cite{Itz}.
Finally,
one deduces
the remaining case $\Sig=P_k$,
$D= dL-\sum_{i=1}^k d_iE_i$,
from the case
$\Sig=\PP^2, D=dL,$ by means of
the inequality
\begin{equation}
N_{P_k,D} \le N_{\PP^2,dL} \left(\prod_{i=1}^kd_k!\right)^{-1}
\quad \text{if} \ c_1(P_k) \cdot D>1 \ \label{enn102}
\end{equation}
applying (\ref{enn102}) to $nD$, $n\ge 2$, instead of $D$, and
using the Styrling formila for the factorial.

We prove (\ref{enn102}) for $P_k$ with any $k\in\N$
assuming, as usual, that
the blown up points are in general position. Namely, we find
a generic configuration $\bz$ of
$$c_1({P_k})\cdot D-1=3d-d_1-...-d_k-1$$ real points in $P_k$ and a
generic configuration $\bz'$
of $3d-1$ real points in $\PP^2$ such that to any
real rational curve in $|D|$
passing through
$\bz$,
correspond $\prod_id_i!$ real plane rational curves of degree $d$
passing through $\bz'$.

Let $\pi:P_k\to\PP^2$ be the
blowing-up at generic real points $p_1,...,p_k\in\PP^2$, and
$\bz$ a collection of $3d-d_1-...-d_k-1$ generic real
points in $P_k$. Since
$c_1(P_k)\cdot D>1$, Theorem 4.1(i) from \cite{GP} applies and it
shows that any element in the set ${\cal R}_D(\bz)$ of rational
curves in $|D|$, passing through $\bz$, is an immersed curve.
Furthermore, slightly moving one of the point of $\bz$ all these
elements become transversal to the exceptional divisors.
Thus, they descend to
plane rational curves of degree $d$ having precisely
$d_i$ nonsingular branches at $p_i$ for any $i=1,\dots, k$.

For each $i=1,\dots, k$, pick a collection $\bz_i$ of $d_i$
real generic points in a
small
neighborhood of $p_i$. For each $C\in {\cal
R}_D(\bz)$, let
$\phi: \PP^1\to \PP^2$ be
an immersion parametrizing
$\pi(C)$. The normal bundle ${\cal N}=\phi^* T\PP^2/T\PP^1$ is a
line bundle of degree $3d-2$ on $\PP^1$, and the linear system
$H^0({\cal N})$ has no base points. Thus, imposing on deformations
of $\phi$ the conditions to pass through all the
$3d-d_1-\dots-d_k-1$ elements of $\bz$, one can freely vary the
images of all the $d_1+\dots+d_k$ points of
$\phi^{-1}\{p_1,\dots,p_k\}\subset \PP^1$,
and therefore we obtain at least $\prod_i d_i !$ real rational curves
passing through $\bz'=\bz\cup\bigcup_i\bz_i$. \proofend

{\small
{\bf Remarks}.

{\bf (1)}
Lemma \ref{upperbound} and its proof extend literally to symplectic
blow-ups of $\PP^2$
as soon as the numbers $N_{nD}$ are replaced by the symplectic
Gromov-Witten
genus zero
invariants.

{\bf (2)} If $c_1(P_k)\cdot D=1$, the inequality (\ref{enn102})
can be replaced by
\begin{equation}
N_{P_k,D}\le N_{\PP^2,(d+1)L}\left(d\prod_{i=1}^kd_k!\right)^{-1}.
\label{dop}\end{equation} The set of rational curves in $|D|$ is
finite, and to prove (\ref{dop}) one can choose as $L$ any line
transversal to all the rational curves $C$ in $|D|$. Pick then
three generic points $z_1,z_2,z_3$ on $L$, and a generic line
$L'\ne L$ passing through $z_3$. A standard application of
Riemann-Roch theorem shows that smoothing up any of the
intersection points of $L$ and $C$, one can include $C+L$ in a
one-dimensional family of rational curves in $|D+L|$ passing
through $z_1,z_2$, and that this family sweeps a neighborhood of
$z_3$ in $L'$. Therefore, there exist at least $dN_{P_k,D}$
rational curves in $|D+L|$ passing through $z_1, z_2$ and a
generic point $z'_3\in L'$ close to $z_3$. Now, it remains to
apply (\ref{enn102}) to $D+L$ instead of $D$.

{\bf (3)} The results of this section provide some information on
convergency properties of the Gromov-Witten potential, which in
the case of $\Sigma=P_k$ we write as the Laurent series
$$
f=\sum_{c_1(P_k)\cdot D\ge 1}\frac{N_D}{(c_1(P_k)\cdot
D-1)!}x^dy_1^{d_1}\cdots y_k^{d_k}, \quad D=dL-\sum d_iE_i.
$$
First of all, the inequalities (\ref{enn102}) and (\ref{dop}) imply that the
convergency domain of $f$ has a nonempty interior. More precisely,
the convergency domain
contains the set
$$
\{|y_1|+\dots +|y_k|<1,\ |x|(1+|y_1|+\dots +|y_k|)^3<a^{-1}\},
$$
where $a=0,138\dots$ is the radius of convergency of $\sum
\frac{N_d}{(3d-1)!}t^d$ (see \cite{Itz}). Indeed, $f$ is term by
term bounded from above, except the finite number of (Laurent)
terms corresponding to $D=E_i$, $1\le i\le k$, by
$$\displaylines{
\sum_{c_1\cdot D\ge 1, (d_1,\dots, d_k)\ge
0}\frac{d^{-1}N_{d+1}}{\prod (d_k)!\cdot(3d-(d_1+\dots
+d_k)-1)!}x^dy_1^{d_1}\cdots y_k^{d_k} \cr
=\sum_{d>0}\frac{d^{-1}N_{d+1}}{(3d-1)!}(1+y_1+\dots
+y_k)^{3d-1}x^d.}
$$

Furthermore, from the lower bound
$$
\log N_{nD}\ge \log W_{nD} \ge (c_1(\Sigma) \cdot D)\cdot n\log n
+ O(n)
$$
it follows that the convergency domain of $f$ is contained in the
intersection of sets
$$
\{|x^dy_1^{d_1}\cdots
y_k^{d_k}|\le R(d,d_1,\dots,d_k)\},
$$
where $R(d,d_1,\dots,d_k)$ are some positive finite
constants, $(d_1,\dots, d_k)\ge 0$ and $c_1(P_k)\cdot D\ge 1$.
This assertion is a routine corollary of the following Cauchy
formula for subdiagonals of a convergent power series $f=\sum
a_{k\hat{l}}x^ky_1^{l_1}\cdots y_k^{l_k}$ in $k+1$ variables $x,
\hat y$ ($\hat y=(y_1,\dots, y_k)$):
$$
\sum_n a_{pn,\hat{q}n}t^{pn}=\frac1{(2\pi
i)^k}\int_{T}f(\frac{t}{\hat y^{\hat q}}, \hat y^p)\frac{d\hat
y}{\hat y}, p\in\N, \hat q\in\N^k,
$$
where
$|t|^p\le x^p(0)\hat y(0)^{\hat q}$,
$T$ is a torus
$\{|y_i|^p=y_i(0)\}$, and $p,\hat q$ are coprime. This integral
formula holds for any positive $(x(0),\hat y(0))$ in the interior
of the convergency domain. The boundedness of $x^p\hat y^{\hat q}$
is thus equivalent to finiteness of the convergency radius of the
subdiagonal, and the latter follows from the above lower bound.}

\subsection{Proof of
the Key Bound (Theorem~\ref{kb}) }\label{proof}

We use the same terminology, notations, and constructions as in
\cite{IKS}. The only difference concerns
the subdivisions called compressing
in \cite{IKS}; here, we prefer to follow the terminology from
general Shustin's scheme \ref{secn10} and call them {\it consistent
subdivisions}, or {\it $\gamma $-consistent subdivisions} when
we want to underline that the subdivisions are obtained by the compressing
procedure applied to a fixed path $\gamma$.

We consider, first, the cases of $\PP^2$ and $Q$, and
then show how to deduce from them the required statements in all
the remaining cases.

\subsubsection{Case $\Sigma = \PP^2$}\label{plane-case}
It is sufficient to treat $D=[\PP^1]$.

\begin{figure}
\begin{center}
\epsfxsize 100mm \epsfbox{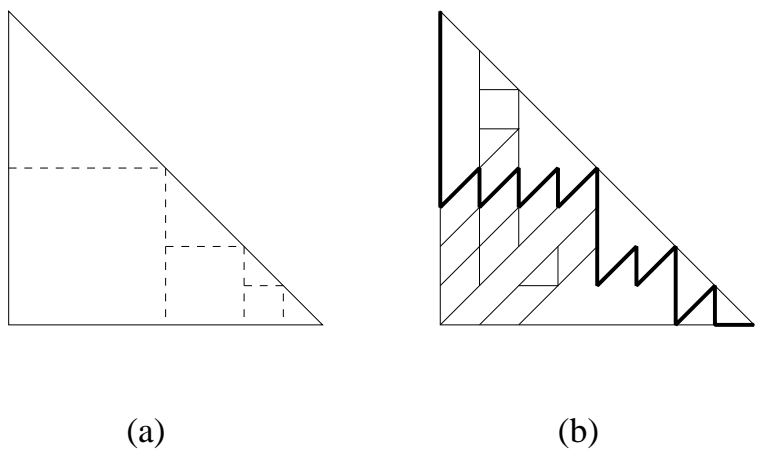}
\end{center}
\caption{Path $\gamma$ and $\gamma$-consistent subdivisions for
$\Sig=\PP^2$} \label{pla}
\end{figure}

Let us fix $n\in\N$ and consider a linear function $\lambda^0:
\R^2 \to \R$ defined by $\lambda^0(i,j) = i - \varepsilon j$,
where $\varepsilon > 0$ is sufficiently small constant (so that
$\lambda^0$ defines a kind of a lexicographical order on the
integer points of the triangle $\Delp=\Delp_{nD}$). Inscribe in
$\Delp$ a sequence of maximal size squares as shown on Figure
\ref{pla}(a). Their right upper vertices have the coordinates
$$
(x_i, y_i),\ i\ge 1,\ x_1=y_1=\left[\frac{n}2\right],\
y_{i+1}=\left[\frac{n-x_i}2\right],\ x_{i+1}=x_i+y_{i+1}.
$$
Put $(x_0,y_0)=(0,n)$. Then pick a $\lambda^0$-admissible lattice
path $\gamma$ consisting of segments of integer length $1$ as
shown on Figure \ref{pla}(b). This path consists of sequences of
vertical segments joining $(x_i, y_i)$ with $(x_i, y_{i+1}-1)$,
zig-zag sequences joining $(x_i, y_{i+1}-1)$ with
$(x_{i+1},y_{i+1})$ (in such a zig-zag sequence the segments of
slope $1$ alternate with vertical segments; it always starts and
ends with segments of slope $1$), and the segments $[(n-1,1),
(n-1,0)]$ and $[(n-1,0), (n,0)]$. The total integer length of this
lattice path is $3n-1$.

Now, we select some $\gamma$-consistent subdivisions of $\Delp$,
that is subdivisions  constructed along the compressing procedure
(see Subsection C in \cite{IKS}) starting with $\gamma$. Namely,
the path~$\gamma$ divides $\Delp$ in two parts: the part
$\Delp_+(\gamma)$ bounded by $\gamma$ and $\partial \Delp_+$ (in
our case, $\partial \Delp_+$ is the hypotenuse of the triangle
$\Delp$), and the part $\Delp_-(\gamma)$ bounded by $\gamma$ and
$\partial \Delp_-$ (two other sides of the triangle). Subdivide
$\Delp_+(\gamma)$ in vertical strips of integer width $1$. Note
that the rightmost strip consists of one primitive triangle. Pack
into each strip but the rightmost one the maximal possible number
of primitive parallelograms and place in the remaining part of the
strip two primitive triangles (see Figure~\ref{pla}(b)). Then
subdivide $\Delp_-(\gamma)$ in slanted strips of slope~$1$ and
horizontal width~$1$. Pack into each strip the maximal possible
number of primitive parallelograms. This gives a subdivision of
any slanted strip situated above the line $y = x$. For any strip
situated below the line $y = x$ place in the remaining part of the
strip one primitive triangle (see Figure~\ref{pla}(b)).

The total number of such
$\gamma$-consistent
subdivisions
is
\begin{equation}\label{compcount}
M_n\ge\prod_i \frac{y_i !(y_i+1)!}{2^{y_i}}\cdot \prod_i {y_i !}\ ,
\end{equation}
where the first product corresponds to subdivisions of
$\Delp_+(\gamma)$ and the second one to those of
$\Delp_-(\gamma)$.

All the constructed subdivisions of $\Delp$ are
nodal and odd (see \cite{IKS} for definitions), each of them is dual to an
irreducible tropical curve and contributes $1$
to the
Welschinger number. For example, the irreducibility of the dual
tropical curve can easily proved by the following induction.
Let us scan the subdivision by vertical lines from right to left.
The rightmost fragment of the tropical curve is dual to the
primitive triangle $\Delp\cap\{x\ge n-1\}$, so it is irreducible.
At the $i$-th step, $i>0$, we look at the irreducible
components of the curve dual to the union of those elements of our
subdivision which intersect the strip $n-i-1<x<n-i$.
Each of these irreducible components either connects the lines
$x=n-i-1$ to $x=n-i$, or contains a pattern dual to a trianlge
with an edge on $x=n-i$, or contains a pattern dual to
a slanted parallelogram. Therefore, each component joins
the curve dual to the subdivision of $\Delp\cap\{x\ge n-i\}$.

According to Proposition~2.5 from \cite{IKS}, any chosen
subdivision contributes $1$ to the Welschinger invariant, and
according to Proposition~2.6 from \cite{IKS} for our choice of the
function $\lambda^0$, the input of any other consistent
subdivision is nonnegative. Therefore, $ W_n\ge M_n $ and it
remains to check that
$$
\log M_n \ge 3n\log n +O(n).
$$

To that purpose, we observe that according to the relations
defining the coordinates $x_i, y_i$ the following inequalities
hold: $y_i\le n-x_i\le 2y_{i+1}+1.$ They imply
$$
y_{i+1}\ge 2^{-i}(y_1+1)-1\ge \frac{n}{2^{i+1}}-1,
$$
where $i$ varies from $0$ to $k\ge [\log_2n]-1$. Then we use the
Styrling relation
$$
\log\Gamma(x)=x\log x+O(x),\
x\to +\infty,
$$
to get finally

$$
\displaylines{
\log \prod {y_i!}\ge
\sum_{i=1}^{[\log_2(n)]}\log\Gamma(
\frac{n}{2^i}-1)
=\sum_{i=1}^{[\log_2(n)]}\frac{n}{2^i}\log\frac{n}{2^i} + O(n)
\cr
= \sum_{i=1}^{\infty}\frac{n}{2^i}\log{n} + O(n)
=n\log n + O(n).
}
$$
\proofend

\subsubsection{Case $\Sigma = Q$}
\label{hyp-case}
\begin{figure}
\begin{center}
\epsfxsize 115mm \epsfbox{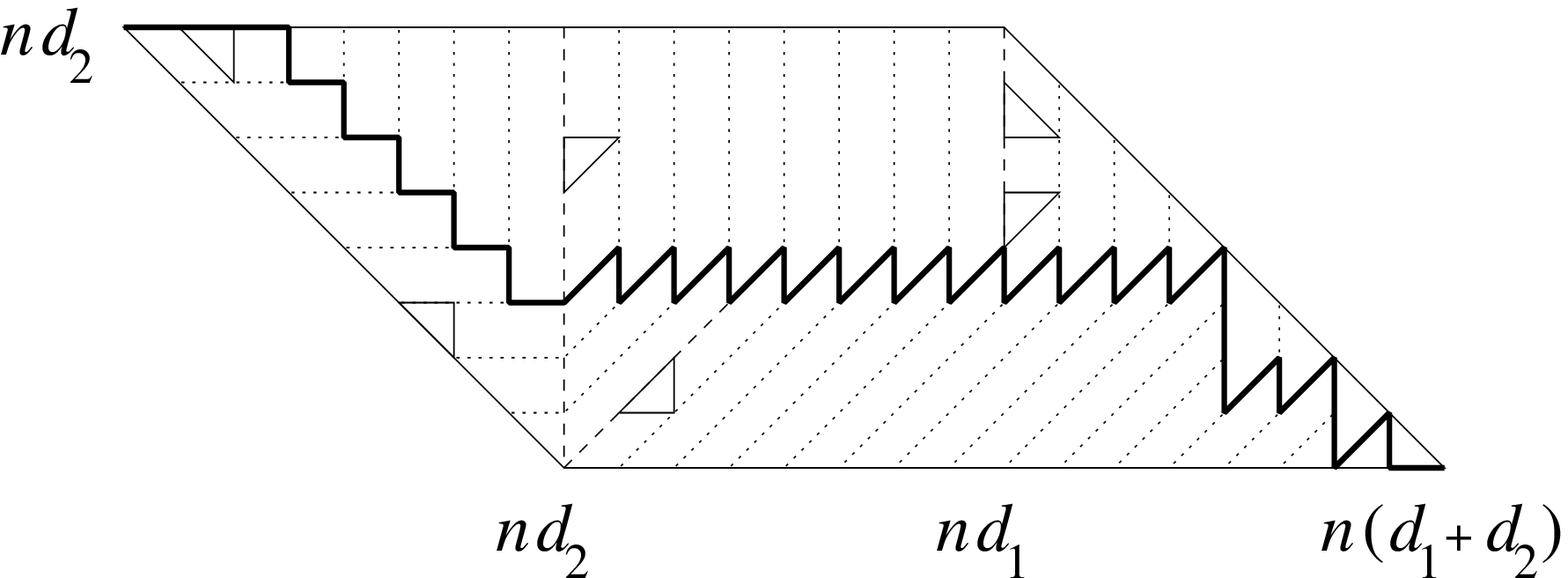}
\end{center}
\caption{Compressing subdivisions for $\Sig=\PP^1 \times \PP^1$}
\label{hyp}
\end{figure}

Assume that $d_1 \geq d_2$. Fix $n$ and consider the
linear function
$\lambda^0: \R^2 \to \R$ introduced in \ref{plane-case}.
Put
$a = [nd_2/2]$. In the
parallelogram
$\Delp=\Delp_{nD}$ pick a $\lambda^0$-admissible lattice path
$\gamma$ consisting of segments of integer length $1$ as shown on
Figure~\ref{hyp}.
This path consists of the sequence of horizontal segments joining
$(0, nd_2)$ with $(a - 1, nd_2)$, a zig-zag sequence joining $(a-
1, nd_2)$ with $(nd_2, a - 1)$
(in this zig-zag sequence the vertical segments
alternate with horizontal ones; it starts with a vertical segment
and ends with a horizontal one),
a zig-zag
sequence joining $(nd_2, a - 1)$ with $(nd_1 + nd_2 - a, a)$
(in this zig-zag sequence the segments of slope $1$
alternate with vertical segments; it  starts and ends with segments
of slope $1$),
and the sequence of segments joining $(nd_1 + nd_2 - a, a)$ with
$(nd_1 + nd_2, 0)$ and obtained by
a translation of the path described in section~\ref{plane-case}
(one should replace $n$ by $a$ in this description).
The total integer
length of this lattice path is $2n(d_1 + d_2) - 1=c_1(Q)\cdot
nD-1$.

Now,
we select
some compressing subdivisions of
$\Delp$ constructed along the compressing procedure starting with
$\gamma$. Namely, subdivide
$\Delp_+(\gamma)$ in vertical strips of width $1$. Pack
into each strip the maximal possible number of primitive
parallelograms and place in the remaining part of the strip one or
two primitive triangles.
The
number of these triangles
is two if the strip is contained in the half-plane
$x \geq nd_1$,
and one
if the strip lies between the vertical lines
$x = nd_2$ and $x = nd_1$; the strip does not contain
any triangle if the strip is contained in the half-plane
$x \leq nd_2$.

The lines $x = nd_2$ and $y = x - nd_2$ divide $\Delp_-(\gamma)$
into three parts (see Figure~\ref{hyp}). Subdivide these parts as
follows:
\begin{itemize}
\item the
left part is subdivided
in
horizontal strips of height~$1$; each strip is
filled in by parallelograms
and one primitive triangle
(see Figure~\ref{hyp});
\item the central part
is subdivided in slanted strips of slope~$1$ and horizontal width~$1$;
each strip is filled in by primitive parallelograms;
\item the right part
is subdivided in similar slanted strips;
each strip is filled in by
primitive parallelograms and one primitive triangle
(see Figure~\ref{hyp}).
\end{itemize}

Since each constructed subdivision contributes $1$ in
the Welschinger invariant and the other consistent subdivisions
have nonnegative contributions, the required bound is obtained
by computations completely similar to those performed in \ref{plane-case}.
\proofend

\subsubsection{Other cases}\label{other-cases}

We start with a reduction of the case $\Sigma = P_1$
to the case
$\Sigma=Q$.
Divide the trapeze
$\Delp=\Delp_{nD}$
by the vertical line $x = n(d - d_1)$
into the rectangle
$\Delp_1$ and the triangle
$\Delp_2$. In the trapeze
$\Delp$ consider $\lambda^0$-admissible lattice paths $\gamma$ of
the form $\gamma = \gamma_1 \cup \gamma_2$, where
\begin{itemize}
\item $\gamma_1$ is any
lattice path in
$\Delp_1$ which starts at the left upper corner of
$\Delp_1$, ends at the right lower corner of
$\Delp_1$, consists of $2nd - 1$ segments of integer length $1$,
and can be completed to a compressing
subdivision of
$\Delp_1$ with a positive input to the Welschinger invariant;
\item $\gamma_2$ is the lattice
path which
consists of $n(d - d_1)$
segments of
length $1$ and goes along the lower
horizontal side of
$\Delp_2$.
\end{itemize}

\begin{figure}
\begin{center}
\epsfxsize 70mm \epsfbox{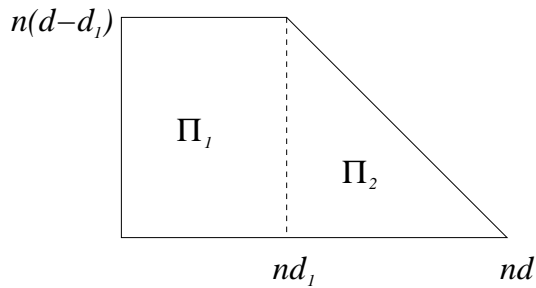}
\end{center}
\caption{Case $\Sig = P_1$} \label{P1}
\end{figure}

Since $\gamma_1$ ends at the right lower corner of $\Delp_1$, any
$\gamma_1$-consistent subdivision of $\Delp_1$ can be completed to
a $\gamma$-consistent subdivision of $\Delp$ in the following way.
Divide $\Delp_2$ into vertical strips of width~$1$. Pack in each
strip the maximal possible number of primitive parallelograms, and
put in the remaining part of each strip one primitive triangle.
This can be done in $(n(d - d_1))!$ ways. Note that if a chosen
$\gamma_1$-consistent subdivision of $\Delp_1$ is dual to an
irreducible tropical curve, then the resulting $\gamma$-consistent
subdivision of $\Delp$ is dual to an irreducible tropical curve as
well. Thus, $W_{nD}$ is bounded from below by $W_{Q,(nd_1,
n(d-d_1))}\cdot(n(d - d_1))!$. It implies the required result in
the case $\Sigma = P_1$.

The case $\Sigma = P_2$
can
be reduced
to the case
$\Sigma=Q$
in a similar way.
The only difference is that the polygon
$\Delp$ is divided into a rectangle and a trapeze.
The case $\Sigma = P_3$ is in its
turn reduced to the case of $\Sig=P_2$. \proofend

\section{Few
results for data
containing imaginary points}\label{further}

In this section we
study the invariants $W_{d,m}$, $m>0$, of toric Del Pezzo
surfaces using as a basis the results of \cite{Sh1}.

\subsection{Examples of calculation
of Welschinger invariant}\label{examples}

(1) All the plane rational cubics interpolating fixed $3d-1=8$ points
belong to the same linear pencil.
The
integration
with respect to the Euler characteristic along the pencil gives
the following {\bf linear}
function
$$W_{3,m}=8-2m, \quad 0\le m\le 4\ .$$

Note, that whatever is $0\le m\le 4$,
the number
$R_{3,m}(\bw)$ of real rational cubics attains the value $W_{3,m}$ for a
suitable generic collection $\bw\in\Omega_{r,m}, r=3d-1=8$.
For example, the pencil spanned by the curves
$x(x^2+y^2+z^2)+\epsilon y^3$ and $y(x^2+y^2+z^2)-\epsilon x^3$,
where $\epsilon$ is sufficiently small,
has $8$ imaginary fixed
points
and does not contain real
singular cubics (to check
this property,
it is sufficient to notice
that
if
$t$ is real, the polynomial $y^3-tx^3$
cannot
vanish at
any
intersection
point
of $x+ty=0$ and $x^2+y^2+z^2=0$).

\noindent (2) To calculate the Welschenger invariant for plane
quartics and quintics, we use birational transformations and
Welschinger's wall crossing formula (see
\cite{W1},
Theorem 2.2)
which expresses
the first finite
difference of the function
$m\mapsto W_{D,m}$
as
twice the Welschinger invariant of the
surface $\Sig $
blown up at
one real point.
For plane quartics  ($\Sigma =\PP^2$, $D=4L$)
we obtain a system of difference
equations
$$
\displaylines{ \Delta^2 W_{4,m}= 4W_{3,m}=4(8-2m), \cr -\Delta^1
W_{4,m}|_{m=0}= 2W_{P^1\times P^1,D=(2,3), m=0}=2\cdot 48=96, \cr
W_{4,m}|_{m=0}=240 }
$$
(two last values are obtained by the algorithm described in
\cite{Mi, Sh0},
see also
\cite{IKS}).
For plane quintics we get similarly
$$
\displaylines{ -\Delta^3 W_{5,m}= 8W_{4,m}, \cr \Delta^2
W_{5,m}|_{m=0}= 4W_{P^1\times P^1, D=(3,3), m=0}=4\cdot 1086, \cr
-\Delta^1 W_{5,m}|_{m=0}= 2W_{P^2(1), D=(5,2), m=0}=9168, \cr
W_{5,m}|_{m=0}=18264. }
$$

Therefore, the Welschenger invariants of plane quartics take the
values
\begin{center}
\begin{tabular}{|l|c|c|c|c|c|c|c|}
\hline $m$ &\makebox[0.35cm] 0 &\makebox[0.35cm] 1
&\makebox[0.35cm] 2&\makebox[0.35cm] 3 &\makebox[0.35cm] 4
&\makebox[0.35cm] 5  \\
\hline
$W$ & 240 & 144 & 80 & 40 & 16 & 0 \\
\hline
\end{tabular}
\end{center}
which are interpolated by a polynomial of degree three
$$
W_{4,m}=-\frac43 m(m-1)(m-2)
+16m(m-1)-96m+240.
$$

For plane quintics the invariants take the values
\begin{center}
\begin{tabular}{|l|c|c|c|c|c|c|c|c|c|}
\hline $m$
&\makebox[0.35cm] 0 &\makebox[0.35cm] 1
&\makebox[0.35cm] 2&\makebox[0.35cm] 3 &\makebox[0.35cm] 4
&\makebox[0.35cm] 5&\makebox[0.35cm] 6 &\makebox[0.35cm] 7 \\
\hline
$W$ & 18264 & 9096 & 4272 & 1872 & 744 & 248 & 64 & 64 \\
\hline
\end{tabular}
\end{center}
which are interpolated by a polynomial of degree six
$$
\displaylines{ W_{5,m}=\frac4{45}m(m-1)(m-2) (m-3)(m-4)(m-5)- \cr
-\frac{32}{15}m(m-1)(m-2) (m-3)(m-4)+ 32m(m-1)(m-2) (m-3)- \cr
-320m(m-1)(m-2)+2172m(m-1) -9168m+18264. }
$$

In these two cases, $d=4$ and $5$,
the degree of the interpolating polynomials
happens
to be
smaller
than
a generic interpolation data,
that is,
smaller than $[\frac{3d-1}2]$.
Let us note that
it is no more the case for any $d\ge 6$.

\noindent (3)
The Welschinger invariants for curves of bi-degree $(d,2)$ on
$\PP^1\times\PP^1$ also can be computed explicitly.
In this case, the number $m$ of pairs of
imaginary
points
varies in the range $0\le m\le d+1$, and
Welschinger's wall
crossing formula takes the form
\begin{equation}
\Del^1 W_{(d+1,2),m}=-2W_{(d,2),m},\quad 0\le m\le d\ .
\label{eh2}\end{equation} A solution of this difference equation
is uniquely determined, for instance, by the initial data
$W_{(1,2),m}=1, 1\le m\le 3$, and the sequence $W_{(d,2),d+1}$,
$d\ge 1$. We claim that
\begin{equation}W_{(d,2),d+1}=\left[\frac{d+1}{2}\right]\cdot 2^{d-1},\quad
d\ge 1\ ,\label{eh1}\end{equation} and prove this formula in
subsection \ref{exh1}.
Thence,
we obtain by induction all the other values:
$$W_{(d,2),m}=(d+m)\cdot 2^{2d-2-m},\quad 0\le
m\le d\ .
$$

For comparison,
$W_{(2,d),0}=d\cdot 2^{2d-2}$ and,
according to \cite{Itz},
$N_{(2,d)}
=d(d+1)\cdot 2^{2d-3}$.

\subsection{
Two general results}\label{further-results}

The
calculations of the Welschinger invariant presented in the previous subsection
lead to the following conjecture.

\begin{conjecture}
Let $\Sig$ be a
toric Del Pezzo surface equipped with its tautological real structure,
and $D$ an ample divisor on $\Sig$.
The Welschinger invariants $W_{D,m}$
are positive if $m<[(c_1(\Sig)\cdot D-1)/2]$,
are non-negative if $m=[(c_1(\Sig)\cdot D-1)/2]$, and satisfy
the monotonicity relation
$$
W_{D,m-1} \ge W_{D,m} \quad\text{for} \quad m\ge 1 \ .
$$
\end{conjecture}

This conjecture is supported by
the following two statements.

\begin{theorem}\label{tn4}
Let $\Sig$ be $\PP^2$, $P_1$, $P_2$, or $P_3$
equipped with its tautological real structure,
$D$ an ample divisor on $\Sigma$, and $m$ is a positive integer less or equal to
$(c_1(\Sig)\cdot D-1)/2$.
The invariants $W_{\Sigma,D,m}$ are positive
for $m \leq 3$ in the case $\Sig = \PP^2$
and for $m \leq 2$ in the cases $\Sig = P_k$, $k = 1, 2, 3$.
Furthermore,
$$\log W_{\Sig, nD, m} =
\log N_{\Sig, nD} + O(n)$$ for $m \leq 3$ in the case $\Sig =
\PP^2$ and for $m \leq 2$ in the cases $\Sig = P_k$, $k = 1, 2,
3$. In addition, if $\Sig$ is $\PP^2$, $P_1$, or $P_2$ and the
projections of $\Delp_D$ on the coordinate axes have at least the
length three, then
$$W_{\Sig, D, 0} > W_{\Sig, D, 1} > W_{\Sig, D, 2}.$$
\end{theorem}

\begin{theorem}\label{tn5}
Let $Q= \PP^1\times\PP^1$ be equipped with its tautological real
structure, and $D$ be a divisor on~$Q$ of bi-degree $(d_1,d_2)$
with positive $d_1,d_2$. Then, for any integer $1 \le m <
d_1+d_2$, the Welschinger invariant $W_{Q,D,m}$ is positive.
Furthermore, let a sequence of integers $m(n)$, $n\ge 1$, satisfy
$$m(n)=\mu n+\psi(n),\quad 0\le m(n)<(d_1+d_2)n,
$$
where
$\mu$ is some real number, and
$\psi(n)$ is a sequence of real numbers such that
$\lim_{n\to\infty}\frac{\psi(n)}{n}=0$.
Then
$$
\log W_{Q, nD,m(n)}\ge (2d_1+2d_2-\mu){n\log n} +O(n+\psi(n)\log
n) )\ .$$
\end{theorem}

\subsection{Proof of Theorems \ref{tn4} and \ref{tn5}}
In~\ref{secn10} and~\ref{signformula}
we denote by $\Delp=\Delp_D$ the convex
lattice polygon introduced in Section \ref{notations} and fix a
non-negative integer $m$ such that $2m\le r=c_1(\Sig)\cdot D-1$.

\subsubsection{Counting scheme
}\label{secn10} Here we shortly describe the way to compute
$W_{D,m}$ worked out in \cite{Sh1} (notation in \cite{Sh1} is
slightly different; there,  $m=r''$ (the number of pairs of
imaginary points in the point data) and $r-2m=r'$ (the number of
real points)).

Pick a linear function $\lam:\R^2\to\R$ which is injective on the
set $\Delp\cap\Z^2$, and choose a set ${\cal I}\subset\{1,2 ,
\ldots , r-m\}$ consisting of $m$ elements. The counting goes
through the following construction of subdivisions of $\Delp$ into
lattice triangles and parallelograms.

Let $p$ and $q$ be the vertices of $\Delp$ such that
$$\lam(p)=\min\lam(\Delp),\quad\lam(q)=\max\lam(\Delp)\ .$$
A lattice path $\gam\subset\Delp$ is called
$(\lam,{\cal I})$\emph{-admissible}, if it has
$r-m+1$
vertices
$v_0,...,v_{r-m}$ such that
$$v_0=p,\quad v_{r-m}=q,\quad\lam(v_i)<\lam(v_{i+1}),\ i=0,...,
r-m-1\ ,$$ and its $r-m$ {\it edges} $\sig_i=[v_{i-1},v_i]$,
$i=1,..., r-m$, satisfy the following restrictions:
\begin{itemize}\item if $i\not\in{\cal I}$, then the
integer
length $|\sig_i|$
of $\sig_i$
is odd,
\item if $i\in{\cal I}$, then either $|\sig_i|$ is even,
or there exists
an integer point $v\in\Delp
\backslash\sig_i$ satisfying the inequalities
$\lam(v_{i-1})<\lam(v)<\lam(v_i)$. \end{itemize}

Pick a $(\lam,{\cal I})$-admissible lattice path $\gam\subset\Delp$.
The following recursive procedure describes the
$\gamma$\emph{-consistent subdivisions} of $\Delp$.

First, for any edge $\sig_i$, $i\in{\cal I}$, of odd integer
length, choose an integer point $v'_i\in\Delp \backslash\sig_i$
with $\lam(v_{i-1})<\lam(v'_i)<\lam(v_i)$ such that the integer
length of $[v_{i-1}, v_i']$ and $[v'_{i}, v_i]$ is odd. Then, put
$$
S_0=\gam\cup\bigcup_iT_i,
$$
where $T_i$ is the triangle with the vertices $v_{i-1},v'_i,v_i$,
and denote by $\tau_0$ the
subdivision of $S_0$ into
segments $\sigma_j$ and triangles $T_i$.
We also put $$W_{\lam,{\cal I}}(\gam,\tau_0)=(-1)^a\prod_i|T_i|,
$$
where $a$ is the number of
integer
points in the interior of $S_0$, and $|T_i|$ stands for
the
\emph{integer}
({\it i.e.}, double Euclidean)
\emph{area} of $T_i$.

After that, we
construct a sequence of contractible sets $S_k\subset\Delp$,
$k\ge 0$, equipped with subdivisions $\tau_k$ into lattice segments,
triangles, and parallelograms, and with numbers
$W_{\lam,{\cal I}}(\gam,\tau_k)$, applying by induction the
following rules:
\begin{itemize}
\item if $S_k$ is convex, we stop the construction,
\item if $S_k$ is not convex, then among
the vertices $v$ of $\tau_k$
such that $v$ is a common
vertex of edges $\sig'$ and $\sig''$ bounding an angle
$<\pi$ outside $S_k$,
we choose the vertex with the minimal value of $\lam$,
\item
then we choose one of the two options: either take the
parallelogram $P \subset \Delp$
generated by $\sig'$ and $\sig''$,
and put
$$
S_{k+1}=S_k\cup P,\quad
W_{\lam,{\cal I}}(\gam,\tau_{k+1})=W_{\lam,{\cal I}}(\gam,\tau_k);
$$
or take the triangle $T$ with the sides $\sig'$ and $\sig''$,
unless $|\sig'|$ and $|\sig''|$ are odd and the length of the
third side of $T$ is even, and put
\begin{equation}\label{eq-new1}
S_{k+1}=S_k\cup T, \quad
W_{\lam,{\cal
I}}(\gam,\tau_{k+1})=W_{\lam,{\cal
I}}(\gam,\tau_k)\cdot(-1)^{a(T)}\cdot A(T),
\end{equation}
where
$$a(T)=\begin{cases}\#(\Int(T)\cap\Z^2)+1,\quad&\text{if}\ |\sig'|
\ \text{and} \
|\sig''|\ \text{are even},\\
\#(\Int(T)\cap\Z^2),\quad&\text{otherwise},\end{cases}$$
$$A(T)=\begin{cases}1,\quad&\text{if}\ |T|\ \text{is odd},\\ |T|\cdot|\sig'|^{-1},\quad&\text{if}
\ |\sig'|\ \text{is even and} \
|\sig''|\ \text{is odd},\\
2|T|\cdot(|\sig'|\cdot|\sig''|)^{-1},\quad&\text{if}\ |\sig'| \
\text{and} \ |\sig''|\ \text{are even},
\end{cases}$$
\item if none of the above options is performable, we stop the
construction.
\end{itemize}
If this (compressing) procedure ends up with $S_k=\Delp$, it
provides a subdivision $\tau_k$ which is convex, {\it i.e.},
consists of the linearity domains of some convex piecewise-linear
function $\nu:\Delp\to\R$. This subdivision is dual to a simple
tropical curve (see \cite{Mi,Sh1}). If the curve is irreducible
and if, in addition, the edges of $\tau_k$ which lie on
$\partial\Delp$ have integer length $\leq 2$, then the subdivision
$\tau_k$ is called a $\gamma$-{\it consistent subdivision} of
$\Delp$. Notice that the sets $S_0,...,S_k$ can uniquely be
restored out of $\cal I$, $\gam$, and $\tau_k$.

The main result of \cite{Sh1}
states that
\begin{equation}\label{eq-new2}
W_{D,m} =\sum_{\gam,\tau}W_{\lam,{\cal I}}(\gam,\tau)\ ,
\end{equation}
{\it where $\gam$ runs over all $(\lam,{\cal I})$-admissible
lattice paths in $\Delp$, and $\tau$ runs over all
$\gamma$-consistent subdivisions of $\Delp$.}

The number $W_{\lam,{\cal I}}(\gam,\tau)$ is called
the {\it Welschinger coefficient}
of $\tau$.

In what follows, the edges of even integer length and the
triangles with all the edges of even integer length are called
\emph{even}. The edges of odd integer length are called
\emph{odd}.

\subsubsection{Even edges expansion and
a sign formula}\label{signformula}

As in Section \ref{proof}, let $\lambda^0: \R^2 \to \R$ be defined
by $\lambda^0(i,j) = i - \varepsilon j$, where $\varepsilon >0$ is
a sufficiently small constant.
Fix a $(\lambda^0, \cal I)$-admissible
lattice path~$\gamma \subset \Delp$.

Note that each
non-vertical side of $\Delp$ is parallel
to a vector $(1, s)$ with $s = 0, \pm 1$.
All lemmas of this subsection are, in fact, valid
for more general polygons, namely for polygons
whose non-vertical sides have integer
\emph{slopes}, that is, are parallel to
vectors $(1, s)$ with $s \in \Z$.

\begin{lemma}\label{AL1}
{\it If $\sig=[(i,j),(i',j')]$, where $i'\ge i$, is an edge of a
$\gamma$-consistent subdivision of $\Delp$, then
\begin{gather}
i'-i\le
\begin{cases}&1,\quad\text{if}\ |\sig|\ \text{is odd},\\
&2,\quad\text{if}\ |\sig|\ \text{is even}.
\end{cases}
\end{gather}\label{AE1}
In particular, in any $\gam$-consistent subdivision
\begin{itemize}
\item
each triangle has a vertical side;
\item
triangles
with at least one odd side
do not contain interior
integer
points;
\item
any
triangle
with a non-vertical even
side
is even.
\end{itemize}}
\end{lemma}

{\bf Proof}. The three statements on triangles follow from the
statement concerning the edges. The latter follows in its turn
from the fact that the integer lengths of the boundary edges of a
$\gamma$-consistent subdivision are at most~$2$. Indeed, if $\sig
\not\subset \partial \Delp$, consider the maximal $k$ such that
$\sig\subset\partial S_k$ (see the construction in \ref{secn10}).
Then $S_{k+1}$ is obtained by adding a triangle or a parallelogram
built on $\sig$ and an adjacent edge. Thus, $\partial S_{k+1}$
contains an edge with the horizontal projection as long as that of
$\sig$ and which, due to the compressing rules, has the same
parity as $\sig$. \proofend

We say that non-vertical edges $\sigma^1$, $\ldots$, $\sigma^l$ of
a $\gam$-consistent subdivision~$\tau$ of~$\Delp$ form a {\it
chain} if, for any $1 \leq j \leq l-1$, the edges $\sigma^j$ and
$\sigma^{j+1}$ are either parallel edges of a parallelogram of
$\tau$ or belong to a triangle of $\tau$ which is not contained
in~$S_0$ (see section~\ref{secn10}). If $\sigma^l$ is contained in
the boundary of~$\Delp$, then the chain $\sigma^1$, $\ldots$,
$\sigma^l$ is called an {\it escaping chain} of $\sigma^1$.

\begin{lemma}\label{new-lemma2}
Let~$\tau$ be a $\gamma$-consistent subdivision of~$\Delp$.
Then
\begin{itemize}
\item any non-vertical edge $\sigma_i$ of~$\gamma$
with $i \not\in {\cal I}$ has exactly two escaping chains;
\item any non-vertical even edge
in~$\gamma$
has exactly two escaping chains;
\item any non-vertical side of any triangle in $S_0$
has exactly one escaping chain;
\item the parity
of edges does not change
along any chain;
\item if $\sigma^1$, $\ldots$, $\sigma^l$
is an escaping chain of $\sigma^1$ such that any edge $\sigma^j$
with $j > 1$ does not belong to $\gamma$, then the slopes of edges
do not increase along the chain;
\item escaping chains of
distinct edges of $S_0$
are disjoint, and any non-vertical edge of~$\tau$
belongs to an escaping chain of some edge of $S_0$.
\end{itemize}
\end{lemma}

{\bf Proof}.
All the statements can easily be derived from
the compressing rules and Lemma
\ref{AL1}.
\proofend

\begin{lemma}\label{lnn1}
Let~$\tau$ be a $\gamma$-consistent subdivision of~$\Delp$. Then
the following sign formula holds
$$
\sign\ W_{\lam^0,{\cal
I}}(\gam,\tau)=(-1)^{\sum s(\sig)},
$$
where $\sig$ runs over all
non-vertical even
edges of $\gam$,
and $s(\sig)$
stands
for the difference of the slopes
of
the terminal edges of two escaping chains of~$\sig$.
\end{lemma}

{\bf Proof}. According to Lemma \ref{AL1}, the triangles of odd
integer area do not contain interior integer points. As to the
number of interior integer points of an even triangle, it is equal
to the half of the integer length of its vertical side diminished
by one, since the length of the horizontal projection of a
non-vertical even side is $2$. Now, the required claim follows
from~(\ref{eq-new1}).\proofend

Denote by $s_{\min}$ the minimum of slopes of non-vertical sides
of~$\Delp$.

\begin{lemma}\label{new-lemma1}
Let~$\tau$ be a $\gamma$-consistent subdivision of~$\Delp$. Then
the slope of any non-vertical edge of~$\tau$ is at least
$s_{\min}$. Moreover, the slope of any non-vertical edge
$\sigma_i$ of~$\gamma$ with $i \not\in {\cal I}$ and of any
non-vertical even edge of~$\gamma$ is at least $s_{\min} + 1$.
\end{lemma}

{\bf Proof}. The both statements are immediate corollaries of
Lemma~\ref{new-lemma2}. Indeed, a non-vertical edge of slope $<
s_{\min}$ in $\tau$ cannot have an escaping chain. If an edge
mentioned in the second statement has the
slope~$s_{\min}$, then the both escaping chains of this edge are
formed by parallel edges. Thus, in this case, the corresponding
tropical curve is reducible. \proofend

\subsubsection{Example}\label{exh1}
As an example of application of the counting scheme given in
\ref{secn10}, we evaluate the Welschinger invariant $W_{(2,d),d+1}
= W_{D,d+1}$ for divisors $D$ of bi-degree $(2,d)$, $d\ge 1$, on
$\Sigma=\PP^1\times\PP^1$.

\begin{figure}
\begin{center}
\epsfxsize 140mm \epsfbox{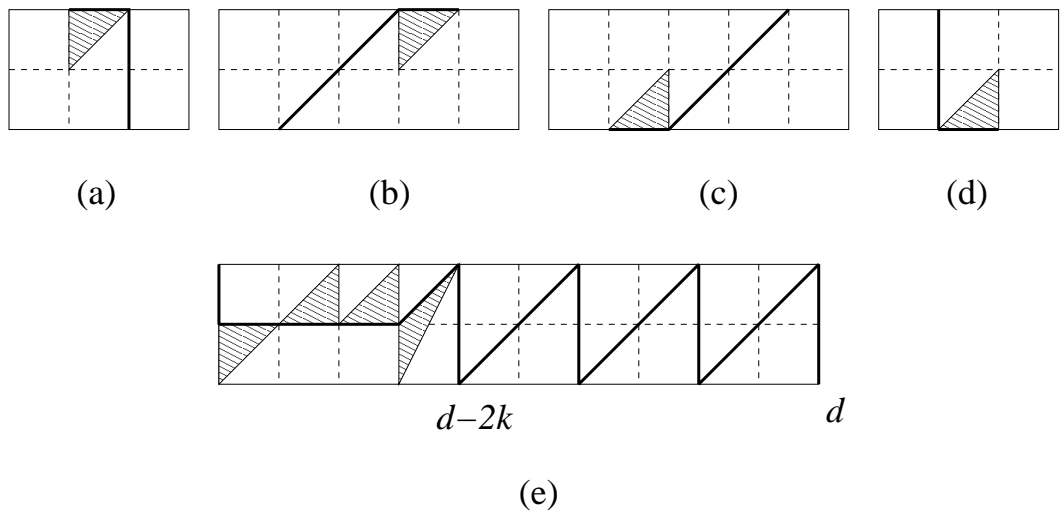}
\end{center}
\caption{Computation of $W_{(2,d),d+1}$} \label{fig29}
\end{figure}

Represent $|D|$ by the rectangle $\Delp$ with vertices $(0,0),
(d,0), (d,2), (2,2)$ and choose ${\cal
I}=\{2,...,d+2\}\subset\{1,...,d+2\}$. Consider a $(\lam^0,{\cal
I})$-admissible lattice path $\gamma$ generating a consistent
subdivision of $\Delp$, and denote the consecutive $d+2$ edges
of~$\gamma$ by $\sig_i$, $1\le i\le d+2$. Since $1 \not\in {\cal
I}$, the edge $\sig_1$ is odd. Then, as it follows from Lemma
\ref{new-lemma1}, the edge $\sig_1$ must be vertical, so of
integer length $1$.
%a segment
Lemma \ref{new-lemma1} implies as well that
each of the other edges
$\sig_i$, $2\le i\le d+2$,
\begin{itemize}
\item
either has integer length $1$,
is of
slope $0$ or $1$, and can be extended to a
triangle of integer area $1$ (as in the construction of $S_0$,
section \ref{secn10}),
\item or has integer length $2$, and is vertical or
of
slope $1$.
\end{itemize}
Therefore, $\gamma$ starts with a sequence of segments of length
$1$ as shown in Figure \ref{fig29}(e); this sequence ends up on
the upper side of $\Delp$. On the other hand, the total integer
length of~$\gamma$ is $d+2$, and its displacement vector is
$(d,-2)$; it implies that $\gamma$ contains exactly one edge of
slope $1$ and integer length $1$, and that the number of vertical
even edges is greater by $1$ than the number of slanted even
edges. In addition, an even edge cannot be followed or preceded by
a horizontal edge of integer length $1$, since otherwise the
subdivision would contain either an odd edge with horizontal
projection of length~$2$, or an edge of negative slope (see
Figures \ref{fig29}(a,b,c,d)). Hence the only remaining
possibilities for $\gamma$ are those shown in Figure
\ref{fig29}(e), where $0\le k\le (d-1)/2$.

Such a path $\gamma$ admits $2^{d-2k-1}$ extensions up to $S_0$:
to any horizontal edge one can attach a triangle either below or
above (see Figure \ref{fig29}(e)). Further extension of $S_0$ to a
$\gamma$-consistent subdivision of $\Delp$ is unique. As a result,
$\Delp\cap\{x\le d-2k\}$ becomes covered by parallelograms,
triangles of integer area $1$, and one triangle of integer area
$2$, while $\Delp\cap\{x\ge d-2k\}$ is covered by triangles of
integer area $4$. The Welschinger coefficient of the subdivision
obtained is $2^{2k}$. Hence,
$$
W_{(2,d),d+1}=\sum_{0\le k\le
(d-1)/2}2^{d-2k-1}\cdot 2^{2k}=\left[\frac{d+1}{2}\right]\cdot
2^{d-1}\ .
$$

\subsubsection{Proof of Theorem \ref{tn4}}
First, we prove the positivity statement.
More precisely, we
show
that, for some specific
choices of ${\cal I}$, the Welschinger coefficients of all
consistent subdivisions of $\Delp$
are positive.

Assume that a
$(\lambda^0, \cal I)$-admissible lattice path $\gamma$ contains at most one
even edge.
Then
any $\gamma$-consistent subdivision
$\tau$
does not contain
any even triangle.
Since, in addition, according to Lemma~\ref{AL1}
the triangles with at least one odd edge
in such a subdivision~$\tau$
do not contain interior integer points, the positivity
of the Welschinger coefficients of all $\gamma$-consistent
subdivisions follows from the formula~(\ref{eq-new1}).
This covers the case $m = 1$.

\begin{figure}
\begin{center}
\epsfxsize 80mm \epsfbox{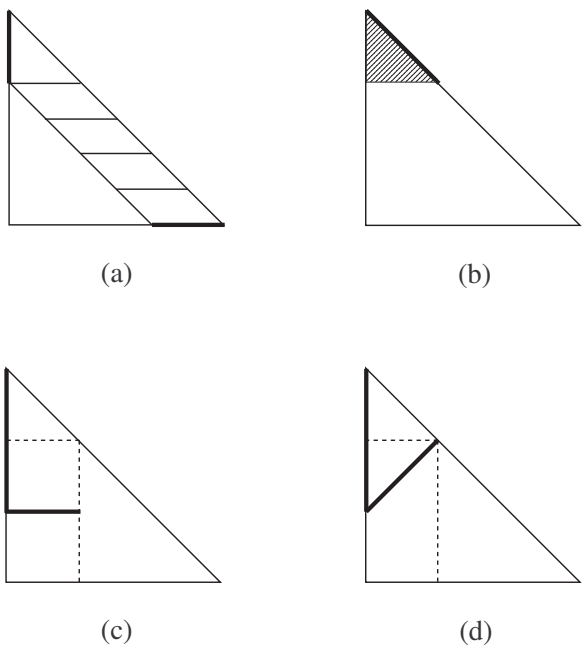}
\end{center}
\caption{Fragments of admissible paths with small $m$}
\label{fign7}

\end{figure}

Let $m = 2$. Put ${\cal I} = \{1, r - m\}$, and assume that a
$(\lambda^0, \cal I)$-admissible lattice path $\gamma$ has exactly
two even edges (the other cases are covered by the previous
consideration). Then the first and the last edges of $\gamma$ are
even. Denote these edges by $l_1$ and $l_2$. Let~$\tau$ be a
$\gamma$-consistent subdivision $\tau$ of $\Delp$ with a negative
Welschinger coefficient. Then $\tau$ has exactly one even
triangle~$T$. Two sides of~$T$ are parallel translations of $l_1$
and $l_2$, and the third side is non-vertical. It implies that
either $l_1$ or $l_2$ is vertical, and thus is contained in the
boundary of~$\Delp$. In addition, the edges of the escaping chains
of the non-vertical edge $l_i$ belong to~$T$ and few
parallelograms (see, for example, Figure \ref{fign7}(a)). Hence,
the tropical curve dual to $\tau$ splits off a component
consisting of three rays with a common vertex. This contradicts to
the consistency of $\tau$ (recall that one of the consistency
conditions is the irreducibility of the dual tropical curve).

Let $m = 3$ and $\Sig=\PP^2$. Put ${\cal I}=\{1,2,3\}$. Consider a
$(\lambda^0, \cal I)$-admissible lattice path~$\gamma$, and assume
that some $\gamma$-consistent subdivision $\tau$ has a negative
Welschinger coefficient. Then, by Lemma \ref{lnn1}, the path
$\gamma$ has either one or three non-vertical even edges. In the
latter case, $\tau$ does not have any even triangle, and according
to Lemma~\ref{AL1} and the formula~(\ref{eq-new1}) the Welschinger
coefficient of~$\tau$ should be nonnegative. Hence, $\gamma$
contains exactly one non-vertical even edge. Assume that the first
edge of~$\gamma$ is odd. As it follows from Lemmas \ref{AL1} and
\ref{new-lemma1}, this edge must be of integer length~$1$ and of
slope~$-1$; see Figure \ref{fign7}(b). The corresponding
triangle~$T$ in $S_0$ has a vertical and a horizontal side.
Lemma~\ref{new-lemma2} implies that all the edges of the escaping
chain of the horizontal side of $T$ are parallel. Thus, in this
case the dual tropical curve splits off a component. Assume now
that the first edge of~$\gamma$ is even. According to
Lemma~\ref{new-lemma1}, this edge is vertical. In addition, it is
of length $2$. If the second edge of~$\gamma$ is odd, then the
same arguments as in the case $m = 2$ show the reducibility of the
tropical curve dual to~$\tau$. Suppose that the second edge
of~$\gamma$ is even. If this edge is non-vertical, it is
horizontal by Lemma~\ref{new-lemma1}, and a component again is
split off the tropical curve; see Figure \ref{fign7}(a). Thus, it
remains to consider the case when the first three edges
of~$\gamma$ are of integer length~$2$, the first and the second
edges are vertical, and the third edge is non-vertical. According
to Lemma~\ref{new-lemma1}, the third edge of~$\gamma$ is either
horizontal or of slope $1$; see Figures \ref{fign7}(c) and
\ref{fign7}(d). In the former case the tropical curve splits off a
component, and in the latter case the third edge cannot generate
an escaping chain in the bottom direction.

\begin{figure}
\begin{center}
\epsfxsize 50mm \epsfbox{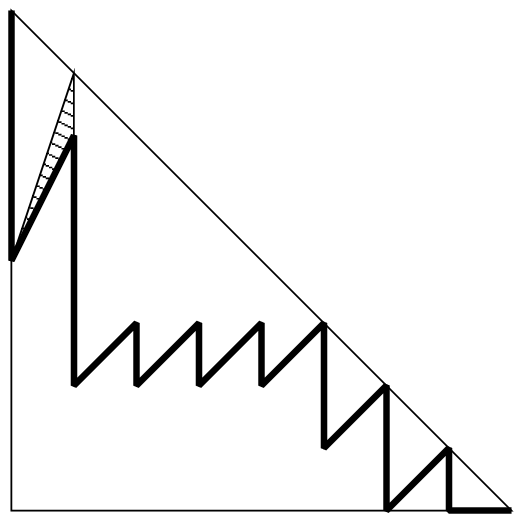}
\end{center}
\caption{Path $\gamma$ for $\Sig=\PP^2$ and $m=2$} \label{plai}
\end{figure}

To prove the asymptotic relations of Theorem \ref{tn4}, we
slightly modify the constructions used in the proof of
Theorem~\ref{kb}. For instance, if $m = 3$ and $\Sig=\PP^2$, we
start the lattice path with two vertical edges $[(0,d),(0,d-2)]$
and $[(0,d-2),(0,d-4)]$ of length~$2$, take the third edge
$[(0,d-4),(1,d-2)]$ of integer length~$1$ supplied with the
triangle with vertices $(0,d-4)$, $(1,d-2)$ and $(1,d-1)$ (see
Figure \ref{plai}), and then continue the path as described in the
proof of Theorem~\ref{kb}.

The final statement of Theorem 4 follows from Welschinger's wall
crossing formula and the positivity of the Welschinger invariant
of $P_1, P_2$, and $P_3$ for the linear system
corresponding to $\Delp_D\setminus \{x+y<2\}$.
\proofend

\subsubsection{Auxiliary bound}
The following lemma is proved exactly
in the same way as Theorem~\ref{kb}.

\begin{lemma}\label{new-lemma3}
The following inequality holds:
$$\log W_{Q, D, 0}\ge 2(d_1+d_2)\log
d_1+O(d_1+d_2),\quad \min\{d_1, d_2\} \to +\infty\ ,
$$
where $D$ is a divisor of bi-degree $(d_1,d_2)$
on $Q=\PP^1\times \PP^1$. \proofend
\end{lemma}

\subsubsection{Proof of Theorem \ref{tn5}}

We prove at the same time the both parts:
the positivity statement
and the lower bound.

Pick an integer $n\ge 1$. In this proof we represent the linear
system $|nD|$ on $Q$ by a rectangle $\Delp$ with vertices $(0,0),
(nd_1,0), (nd_1,nd_2), (0,nd_2)$.
According to Lemma~\ref{lnn1},
for any
$m(n)$-element set ${\cal I}\subset\{1,\dots,2(nd_1+nd_2)-1-m(n)\}$ and
any ($\lambda^0,\cal I)$-admissible path~$\gamma$ in~$\Delp$,
all the $\gamma$-consistent subdivisions of~$\Delp$ have positive
Welschinger coefficients.
Thus,
to prove the theorem, it is sufficient to
construct
an appropriate number
of consistent subdivisions.

Our construction ramifies in few cases.
We use the fact that, as it follows from Lemma \ref{new-lemma1}, if
$1\notin\cal I$
then the first edge of any ($\lambda^0,\cal I)$-admissible path
producing a consistent subdivision is vertical.

Suppose, first, that $m(n)=nd_1+nd_2-1$,
and put ${\cal I}=\{2,\dots,m(n)+1\}$.

\begin{figure}
\begin{center}
\epsfxsize 100mm \epsfbox{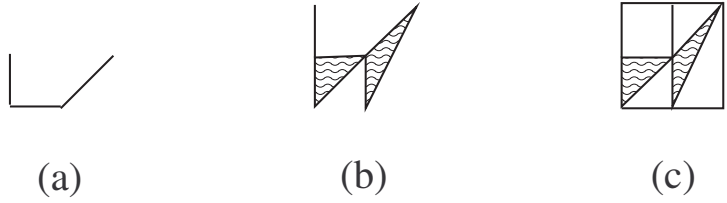}
\end{center}
\caption{Case of $nd_1$ and $nd_2$ even} \label{fig-n1}
\end{figure}

If
$nd_1$ and $nd_2$ are even,
we
consider the paths~$\gamma$ of the form
$\gamma=\gamma_1\cup\gamma_2$, where $\gamma_1$
consists of three edges of integer length~$1$ and
looks like in Figure \ref{fig-n1}(a),
and $\gamma_2$ is obtained by multiplication by $2$ (and translation)
of a $(\lambda^0,\emptyset)$-admissible path $\gamma'_2$ related
to divisors of bi-degree $(nd_1/2-1, nd_2/2)$ and
providing a consistent subdivision.
To obtain the set $S_0$
we complete the second and the third edges of $\gamma_1$
by primitive triangles as is shown in Figure \ref{fig-n1}(b).
Then, we complete $S_0$ to $\gamma$-consistent subdivisions
of $\Delp$ which start from a pattern shown
in Figure~\ref{fig-n1}(c).
Namely, the vertical strip of width~$2$
below the described square
is filled up by rectangles
$1 \times 2$ ({\it i.e.},
of width $1$ and height $2$),
and the remaining rectangle
$(nd_1-2) \times nd_2$
is filled up
by the doubles of
the elements of $\gamma'_2$-consistent subdivisions.

\begin{figure}
\begin{center}
\epsfxsize 100mm \epsfbox{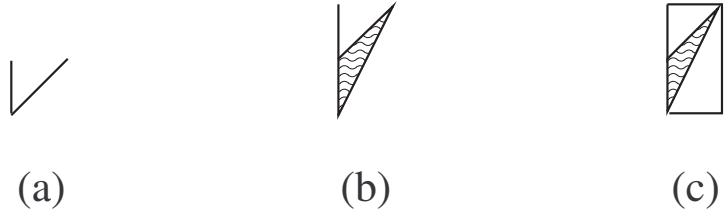}
\end{center}
\caption{Case of $nd_1$ odd and $nd_2$ even} \label{fig-n2}
\end{figure}

If $nd_1$ is odd and  $nd_2$ is even, we consider the
paths~$\gamma$ of the form
$\gamma=\gamma_1\cup\gamma_2$, where $\gamma_1$
consists of two edges of integer length~$1$ and
looks like in Figure \ref{fig-n2}(a),
and $\gamma_2$ is obtained by multiplication by $2$ (and translation)
of a $(\lambda^0,\emptyset)$-admissible path $\gamma'_2$ related
to divisors of bi-degree $((nd_1-1)/2, nd_2/2)$ and
providing a consistent subdivision.
To obtain the set $S_0$
we complete the second edge of $\gamma_1$
by a primitive triangle as is shown in Figure \ref{fig-n1}(b).
Then, we complete $S_0$ to $\gamma$-consistent subdivisions
of $\Delp$ which start from a pattern
shown in Figure~\ref{fig-n2}(c).
Namely, the vertical strip of width~$1$
below the described pattern
is filled up by rectangles
$1 \times 2$,
and the remaining rectangle
$(nd_1-1) \times nd_2$
is filled up
by the doubles of
the elements of $\gamma'_2$-consistent subdivisions.

\begin{figure}
\begin{center}
\epsfxsize 100mm \epsfbox{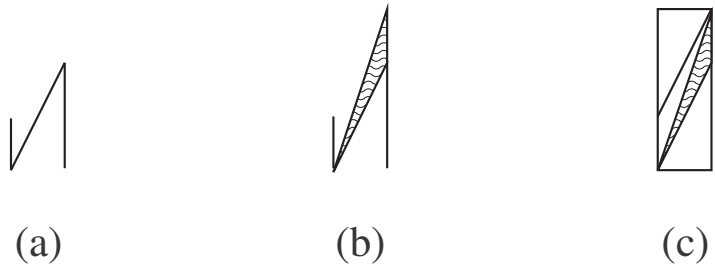}
\end{center}
\caption{Case of $nd_1$ and $nd_2$ odd} \label{fig-n3}
\end{figure}

If
$nd_1$ and $nd_2$ are odd, we
put ${\cal I}=\{1,\dots,m(n)-2,m(n), m(n)+1\}$ and consider
the paths~$\gamma$ of the form
$\gamma=\gamma_1\cup\gamma_2$, where $\gamma_2$
consists of three edges (two edges of integer length~$1$
and one edge of integer length~$2$)
and looks like in Figure \ref{fig-n3}(a),
and $\gamma_1$ is obtained by multiplication by $2$ (and translation) of a
$(\lambda^0,\emptyset)$-admissible path~$\gamma'_1$
related
to divisors of bi-degree $((nd_1-1)/2, (nd_2-1)/2)$ and
providing a consistent subdivision.
To obtain the set $S_0$
we complete the second edge of $\gamma_2$
by a primitive triangle as is shown in Figure \ref{fig-n3}(b).
Then, we complete $S_0$ to $\gamma$-consistent subdivisions
of $\Delp$ in the following way:
we put the doubles of
the elements of $\gamma'_1$-consistent subdivisions
in the rectangle with vertices
$(0,1)$, $(nd_1,1)$, $(nd_1, nd_2)$, and $(0,nd_2)$;
continue with the pattern shown in Figure~\ref{fig-n3}(c);
and fill up the strip situated to the left (resp., above)
of the pattern by rectangles
$2 \times 1$ (resp., $1 \times 2$).

Suppose now that $m(n)<nd_1+nd_2-1$,
and put ${\cal I}=\{j, j+2, j+3, \dots,
j + m(n)\}$, where $j = 2(nd_1 + nd_2 - m(n)) - 1$.

\begin{figure}
\begin{center}
\epsfxsize 100mm \epsfbox{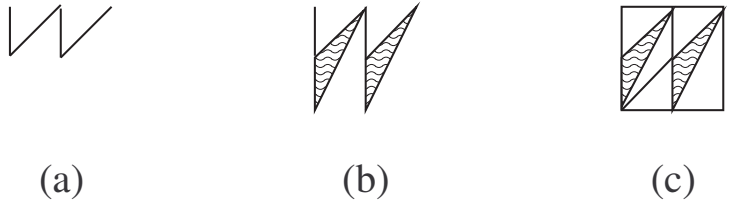}
\end{center}
\caption{Case of odd $m$}
\label{fig-n4}
\end{figure}

If $m(n)$ is odd,
we consider the paths~$\gamma$ of the form
$\gamma=\gamma_1\cup\gamma_2\cup\gamma_3$, where
\begin{itemize}
\item $\gamma_1$ is the translation of a
$(\lambda^0,\emptyset)$-admissible path
$\gamma'_1$ providing a consistent subdivision and related
to divisors of nonnegative bi-degree $(nd_1-2b_1,nd_2-2b_2)$
with $b_1$ and $b_2$ nonnegative integers such that
$2b_1+2b_2=m(n)+1$,
\item $\gamma_2$
consists of four edges of integer length~$1$
and
looks like in Figure \ref{fig-n4}(a),
\item $\gamma_3$ is obtained by multiplication by $2$ (and translation)
of a $(\lambda^0,\emptyset)$-admissible path $\gamma'_3$
providing a consistent subdivision and
related to divisors of bi-degree $(b_1-1, b_2)$.
\end{itemize}
To obtain the set $S_0$
we complete
the second and the forth edges of $\gamma_2$ by
primitive triangles as
in Figure \ref{fig-n4}(b).
Then, we complete $S_0$ to $\gamma$-consistent subdivisions
of $\Delp$ in the following way:
we shift the elements of $\gamma'_1$-consistent subdivisions
in the rectangle with vertices
$(0,2b_2)$, $(nd_1-2b_1,2b_2)$, $(nd_1-2b_1,nd_2)$, and $(0,nd_2)$;
superpose over~$\gamma_2$ the pattern shown on Figure~\ref{fig-n3}(c);
put the doubles of
the elements of $\gamma'_3$-consistent subdivisions
in the rectangle with vertices
$(nd_1-2b_1+2, 0)$, $(nd_1, 0)$, $(nd_1, 2b_2)$,
and $(nd_1-2b_1+2, 2b_2)$;
and fill up the remaining parts by rectangles
of sizes $1 \times 1$, $1 \times 2$, and $2 \times 1$.

\begin{figure}
\begin{center}
\epsfxsize 100mm \epsfbox{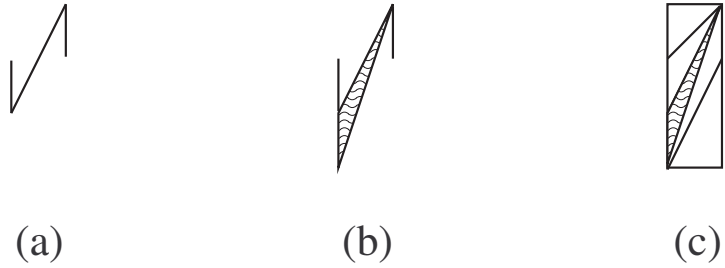}
\end{center}
\caption{Case of even $m(n)$} \label{fig-n5}
\end{figure}

If $m(n)$ is even,
we consider the paths~$\gamma$ of the form
$\gamma=\gamma_1\cup\gamma_2\cup\gamma_3$, where
\begin{itemize}
\item $\gamma_1$ is the translation of a
$(\lambda^0,\emptyset)$-admissible path
$\gamma'_1$ providing a consistent subdivision and related
to divisors of nonnegative bi-degree $(nd_1-2b_1-1,nd_2-2b_2)$
with $b_1$ and $b_2$ nonnegative integers such that
$2b_1+2b_2=m(n)$,
\item $\gamma_2$
consists of three edges of integer length~$1$
and
looks like in Figure \ref{fig-n5}(a),
\item $\gamma_3$ is obtained by multiplication by $2$ (and translation)
of a $(\lambda^0,\emptyset)$-admissible path $\gamma'_3$
providing a consistent subdivision and
related to divisors of bi-degree $(b_1, b_2)$.
\end{itemize}
To obtain the set $S_0$
we complete
the second edge of $\gamma_2$ by
a primitive triangle as
in Figure \ref{fig-n5}(b).
The further construction of
$\gamma$-consistent subdivisions of~$\Delp$
uses the pattern shown on Figure~\ref{fig-n5}(c)
and is performed in the same way as in the case
of odd $m(n)$.

The existence of $\gamma$-consistent subdivisions implies the
positivity of the Welschinger invariants. To prove the required
asymptotic inequality we bound from below, using
Lemma~\ref{new-lemma3}, the number of constructed subdivisions.
For example, in the very last case we take $b_i$, $i=1,2$, close
to $nd_i\cdot\frac{m(n)}{2n(d_1+d_2)}$ and verifying the
conditions imposed in the construction, and then deduce from
Lemma~\ref{new-lemma3} that
$$
\log W_{(nd_1-2b_1-1, nd_2-2b_2),0}\ge (2nd_1+ 2nd_2-2\mu n)\log n
+ O(n+\psi(n)\log n),
$$
if $\mu<d_1+d_2$, and
$$
\log W_{(b_1,b_2),0}\ge \mu n \log n +O(n+\psi(n)\log n),
$$
if $\mu>0$.
Using the positivity of the
Welschinger
coefficients we obtain
$$
\displaylines{
\log W_{D,m(n)}\ge \log W_{(nd_1-2b_1-1, nd_2-2b_2),0}
+ \log W_{(b_1,b_2),0} \cr
\ge
(2d_1+2d_2-\mu) n \log n +O(n+\psi(n)\log n) .}
$$
\proofend

{\ncsc Universit\'e Louis Pasteur et IRMA \\[-21pt]

7, rue Ren\'e Descartes, 67084 Strasbourg Cedex, France} \\[-21pt]

{\it E-mail address}: {\ntt itenberg@math.u-strasbg.fr}

\vskip10pt

{\ncsc Universit\'e Louis Pasteur et IRMA \\[-21pt]

7, rue Ren\'e Descartes, 67084 Strasbourg Cedex, France} \\[-21pt]

{\it E-mail address}: {\ntt kharlam@math.u-strasbg.fr}

\vskip10pt

{\ncsc School of Mathematical Sciences \\[-21pt]

Raymond and Beverly Sackler Faculty of Exact Sciences\\[-21pt]

Tel Aviv University \\[-21pt]

Ramat Aviv, 69978 Tel Aviv, Israel} \\[-21pt]

{\it E-mail address}: {\ntt shustin@post.tau.ac.il}


\begin{thebibliography}{99}

\bibitem{Itz}
Di Francesco, P., Itzykson, C.:
Quantum intersection rings.
{\it in ``The moduli space of curves''}
(Texel Island, 1994),  81--148, Progr. Math., 129,
Birkh\"auser Boston, Boston, MA, 1995.



\bibitem{GP} G\"ottsche, L, and Pandharipande, A.:
The quantum cohomology of
blow-ups of $\PP^2$ and enumerative geometry. {\it J. Diff. Geom.}
{\bf 48} (1998), no. 1, 61--90.

\bibitem{IKS} Itenberg, I., Kharlamov, V., and Shustin, E.: Welschinger
invariant and enumeration of real rational curves. {\it
International Math. Research Notices} {\bf 49} (2003), 2639--2653.

\bibitem{KM} Kontsevich, M.
and Manin, Yu.: Gromov-Witten classes, quantum cohomology and
enumerative geometry. {\it Commun. Math. Phys.} {\bf 164} (1994),
525-562.

\bibitem{Mi0} Mikhalkin,~G.: Counting curves via the lattice paths
in polygons. {\it C. R. Acad. Sci. Paris, S\'er. I}, {\bf 336}
(2003), no. 8, 629--634.

\bibitem{Mi} Mikhalkin,~G.: {\it  Enumerative tropical algebraic geometry in
$\R^2$}. Preprint arXiv: math.AG/0312530.

\bibitem{Sh0} Shustin,~E.: {\it Patchworking singular algebraic
curves, non-Archimedean amoebas and enumerative geometry}.
Preprint arXiv: math.AG/0211278.

\bibitem{Sh1} Shustin,~E.: {\it A tropical calculation of the Welschinger invariants
of real toric Del Pezzo surfaces}. Preprint arXiv:math.AG/0406099.

\bibitem{Vi} Viro,~O.: {\it Dequantization of Real Algebraic Geometry on a Logarithmic Paper}.
Proceedings of the 3rd European Congress of Mathematicians, Birkh\"auser, Progress in Math, 201, (2001), 135--146.

\bibitem{W} Welschinger,~J.-Y.: Invariants of real rational symplectic
4-manifolds and lower bounds in real enumerative geometry. {\it C.
R. Acad. Sci. Paris, S\'er. I,} {\bf 336} (2003), 341--344.

\bibitem{W1} Welschinger,~J.-Y.:
{\it Invariants of real symplectic
4-manifolds and lower bounds in real enumerative geometry},
Preprint arXiv: math.AG/0303145.

\end{thebibliography}
\end{document}